\input ./style/arxiv-general.cfg
\documentclass[MSNbibl,number,citesort,seceqn,dvips]{arxbj}
\makeatletter
   \@ifpackageloaded{graphicx}{}{\usepackage{graphicx}}
\makeatother
\usepackage{mathbh}

\volume{22}
\issue{2}
\pubyear{2016}
\firstpage{969}
\lastpage{994}
\doi{10.3150/14-BEJ683}
\docsubty{FLA}

\makeatletter
\newcommand{\rrvert}{\vert}
\newcommand{\rrVert}{\Vert}
\newcommand{\llvert}{\vert}
\newcommand{\llVert}{\Vert}
\renewcommand{\mid}{|}
\newproclaim{de}{Definition}
\newtheorem{theo}{Theorem}
\newtheorem{prop}[theo]{Proposition} 
\newtheorem{lemma}{Lemma}
\newremark{rem}{Remark}
\makeatother

\begin{document}
\begin{frontmatter}

\title{Sharp ellipticity conditions for ballistic behavior of random
walks in random~environment}
\runtitle{Sharp ellipticity conditions for ballistic behavior}

\begin{aug}
\author[A]{\inits{\'{E}.}\fnms{\'{E}lodie}~\snm{Bouchet}\corref{}\thanksref{A}\ead[label=e1]{bouchet@math.univ-lyon1.fr}},
\author[B]{\inits{A.F.}\fnms{Alejandro F.}~\snm{Ram\'irez}\thanksref{B}\ead[label=e2]{aramirez@mat.puc.cl}}
\and
\author[A]{\inits{C.}\fnms{Christophe}~\snm{Sabot}\thanksref{A}\ead[label=e3]{sabot@math.univ-lyon1.fr}}
\address[A]{Universit\'e de Lyon, Universit\'e Lyon 1, Institut
Camille Jordan, CNRS UMR 5208, 43, Boulevard du 11 novembre 1918, 69622
Villeurbanne Cedex, France.\\
\printead{e1,e3}}
\address[B]{Facultad de Matem\'{a}ticas, Pontificia Universidad Cat\'
{o}lica de Chile, Vicu\~{n}a Mackenna 4860, Macul, Santiago, Chile.
\printead{e2}}
\end{aug}

%
\received{\smonth{2} \syear{2014}}
%
\revised{\smonth{8} \syear{2014}}

%
\begin{abstract}
We sharpen ellipticity criteria for random walks in i.i.d. random
environments introduced by
Campos and Ram\'\i rez which ensure ballistic behavior. Furthermore, we
construct new examples
of random environments for which the walk satisfies the polynomial
ballisticity criteria of
Berger, Drewitz and Ram\'\i rez. As a corollary, we can exhibit a new
range of values
for the parameters of Dirichlet random environments in dimension $d=2$
under which the
corresponding random walk is ballistic.
\end{abstract}

%
\begin{keyword}
\kwd{Dirichlet distribution}
\kwd{max-flow min-cut theorem}
\kwd{random walk in random environment}
\kwd{reinforced random walks}
\end{keyword}
\end{frontmatter}

\section{Introduction}

We continue the study initiated in \cite{CR} sharpening the ellipticity
criteria which ensure ballistic behavior of random walks in random
environment. Furthermore, we apply our results to exhibit a new
class of ballistic random walks in Dirichlet random environments
in dimension $d=2$.

For $x \in\mathbb{R}^d$, denote by $\llvert x\rrvert _1$ and
$\llvert x\rrvert _2$ its $L^1$ and
$L^2$ norm, respectively. Call $U:=\{ e \in\mathbb{Z}^d\dvt \llvert
e\rrvert _1 = 1 \}
= \{ e_1, \ldots,e_{2d} \}$ the canonical vectors with the convention
that $ e_{d+i} = -e_i $ for $1\leq i\leq d$. We set ${\mathcal P}:=\{
p(e)\dvt p(e)\geq0,\sum_{e\in U}p(e)=1\}$.

An environment is an element $\omega:= \{ \omega(x)\dvt x \in\mathbb
{Z}^d \}$ of the environment space $ \Omega:= \mathcal{P} ^{\mathbb
{Z}^d}$. We denote the components of $\omega(x)$ by $\omega(x,e)$.

The random walk in the environment $\omega$ starting from $x$ is the
Markov chain $ \{ X_n\dvt n \geq0 \} $ in $\mathbb{Z}^d$ with law
$P_{x,\omega}$ defined by the condition $ P_{x,\omega}( X_0 = x )=1 $
and the
transition probabilities
\[
P_{x,\omega}(X_{n+1} = x+e \mid X_n = x) = \omega(x,e)
\]
for each $x \in\mathbb{Z}^d$ and $e\in U$.

Let $\mathbb{P}$ be a probability measure defined on the environment
space $\Omega$ endowed with its Borel \mbox{$\sigma$-}algebra, such that $
\{ \omega(x)\dvt x \in\mathbb{Z}^d \} $ is i.i.d. under $\mathbb{P}$.
We call $P_{x,\omega}$ the quenched law of the random walk in random
environment (RWRE) starting from $x$, and $ P_x:= \int P_{x,\omega}
\,\mathrm{d}\mathbb{P}(\omega) $ the averaged or annealed law of the RWRE
starting from $x$.

The law $\mathbb{P}$ is said to be elliptic if for every $ x \in
\mathbb{Z}^d $ and $ e \in U $, $ \mathbb{P}( \omega(x,e) > 0 ) = 1
$. We say that $\mathbb{P}$ is uniformly elliptic if there exists a
constant $\gamma>0$ such that for every $ x \in\mathbb{Z}^d $ and $
e \in U $, $ \mathbb{P}( \omega(x,e) \geq\gamma) = 1 $.

Given $ l \in\mathbb{S} ^{d-1}$ we say that the RWRE is transient in
direction $l$ if
\[
P_0 (A_l) = 1,
\]
with
\[
A_l:= \Bigl\{ \lim_{n \to\infty} X_n \cdot l
= \infty\Bigr\}.
\]
Furthermore, it is ballistic in direction $l$ if $P_0$-a.s.
\[
\liminf_{ n \to\infty} \frac{ X_n \cdot l }{ n } > 0.
\]

Given $ \Lambda\subset\mathbb{Z}^d$, we denote its outer boundary by
\[
\partial\Lambda:= \bigl\{ x \notin\Lambda\dvt \llvert x-y\rrvert_1 = 1
\mbox{ for some } y \in\Lambda\bigr\}.
\]

We denote any nearest neighbour path with $n$ steps joining two points
$ x,y \in\mathbb{Z}^d $ by $ ( x_1, x_2,\ldots, x_n ) $, where
$x_1=x$ and $x_n=y$.

\subsection{Polynomial condition, ellipticity condition}

In \cite{BDR}, Berger, Drewitz and Ram\'{\i}rez introduced a
polynomial ballisticity condition within the uniformly elliptic
context, which was later extended to the elliptic case by Campos and
Ram\'\i rez in \cite{CR}. This condition will be of interest for our
results. It is effective, in the sense that it can a priori be verified
explicitly for a given environment.

To define it, we need for each $ L,\tilde L > 0 $ and $ l \in\mathbb S
^{d-1} $ to consider the box
\[
B _{l,L,\tilde L}:= R \bigl( ( -L,L ) \times( -\tilde L,\tilde L )
^{d-1} \bigr) \cap\mathbb{Z}^d,
\]
where $R$ is a rotation of $\mathbb{R}^d$ that verifies $ R(e_1) = l $.

For each subset $A\subset\mathbb Z^d$ we denote the first exit time
from the set $A$ as
\[
T_A:= \min\{ n \geq0\dvt X_n \notin A \}.
\]
Define also the half space
\[
H_{l}:= \bigl\{ l' \in\mathbb{R}^{d}\dvt
l' \cdot l \geq0 \bigr\}.
\]
We now choose $\alpha>0$ such that
\[
\eta_\alpha:= \max_{e \in H_l\cap\mathbb Z^d} \mathbb{E} \biggl(
\frac{1}{ \omega(0,e) ^\alpha
} \biggr) <\infty
\]
and let
\begin{equation}
\label{ce0} c_0:= \tfrac{2}{3} 3 ^{ 120 d^4 + 3000 d
(({1}/{\alpha}) \log\eta_{ \alpha} ) ^2 }.
\end{equation}

\begin{de}
Given $ M \geq1 $, we say that condition $ (P) _M $ in direction $l$
is satisfied (also written as $ (P) _M \mid l $) if there exists $ L
\geq
c_0 $ and $\tilde L \leq70 L^3 $ such that one has the following upper
bound for the probability that the walk does not exit the box
$B_{l,L,\tilde L}$ through its front side:
\[
P _0 ( X_{ T_{ B_{l,L,\tilde L} } } \cdot l < L ) \leq\frac{1}{L^M}.
\]
\end{de}

This condition has proven to be useful in the uniformly elliptic case.
Indeed, $(P)_M$ for $ M \geq15 d+5 $ implies ballisticity (see \cite{BDR}).

For non-uniformly elliptic environments in dimensions $d\geq2$, there
exist elliptic random walks which are transient in a given direction
but not ballistic in that direction (see, e.g., Sabot--Tournier
\cite{ST}, Bouchet \cite{Bouchet}).
In \cite{CR}, Campos and Ram\'{\i}rez introduced ellipticity criteria
on the law of the environment
which ensure ballisticity if condition $(P)_M$ is satisfied for $M\geq15d+5$
In this article we will sharpen this ellipticity criteria.

\begin{rem}
In definition (1.6) of \cite{CR}, an incorrect value of the constant
$c_0$ is given, different from the definition in (\ref{ce0}).
Nevertheless, it is straightforward to check that the argument of
Section~3.1 of \cite{CR} showing that $(P)_M$ implies $(T)_{\gamma
_L}$ does not change.
\end{rem}

Let us first recall the ellipticity condition of \cite{CR}. For all $M
\geq1$, the polynomial condition $(P)_M$ implies the existence of an
asymptotic direction (see, e.g., Simenhaus \cite{Si}): there
exists $ \hat v \in\mathbb{S} ^{d-1} $ such that $P_0$-a.s.,
\[
\lim_{n\to\infty} \frac{X_n}{\llvert X_n\rrvert _2} = \hat v.
\]
We call $ \hat v $ the asymptotic direction.

\begin{de}
Let $ \beta> 0 $. We say that the law of the environment satisfies the
ellipticity condition $ (E') _\beta$ if there exists an $ \{ \alpha
(e)\dvt e \in U \} \in(0,\infty)^{2d} $ such that
\begin{equation}
\label{kappa} \kappa\bigl( \bigl\{ \alpha(e)\dvt e \in U \bigr\} \bigr
):= 2 \sum
_{e'} \alpha\bigl(e'\bigr) - \max
_{e\in U} \bigl( \alpha(e) + \alpha(-e) \bigr) > \beta
\end{equation}
and for every $ e \in U $
\begin{equation}
\label{finiteness} \mathbb{E} \biggl( \prod_{e}
\omega(0,e) ^{-\alpha(e)} \biggr) < \infty.
\end{equation}
Furthermore, when $ \hat v $ exists, we say that the ellipticity
condition $ (E') _\beta$ is satisfied towards the asymptotic direction
if there exists an $\{ \alpha(e)\dvt e \in U \}$ satisfying (\ref
{kappa}) and (\ref{finiteness}) and such that there exists $ \alpha_1
> 0 $ that satisfies $ \alpha(e) = \alpha_1 $ for $ e \in H _{ \hat
v} \cap U $ while $ \alpha(e) \leq\alpha_1 $ for $ e \in U \setminus
H _{ \hat v} $.
\end{de}

\begin{rem}
In \cite{CR}, (\ref{finiteness}) is replaced by $ \mathbb{E} (
\prod_{e'\neq e} \omega(0,e') ^{-\alpha(e')} ) < \infty$.
Those two conditions are in fact equivalent. The direct implication is
straightforward. And since $ 1 \leq\sum_{e \in U} \mathbh{1} _{ \{
\omega(0,e) \geq {1}/(2d) \} } $, we get
\[
\mathbb{E} \biggl( \prod_{e} \omega(0,e)
^{-\alpha(e)} \biggr) \leq\sum_{e \in U}
(2d)^{ \alpha(e) } \mathbb{E} \biggl( \prod_{e'\neq
e}
\omega\bigl(0,e'\bigr) ^{-\alpha(e')} \biggr) < \infty.
\]
This gives the reverse implication.
\end{rem}

\begin{rem}
\label{eta-alpha}
Knowing the existence of $ \hat v $ does not mean that we know its
value. In most cases, $ \hat v $~is found to be inaccessible. A notable
exception is the result of Tournier \cite{T2} that gives the value
of $ \hat v $ in the case of random walks in Dirichlet environments.
\end{rem}

\subsection{Ballisticity results}

Our main results are a generalization of Theorems 1.2 and 1.3 of \cite
{CR} where we remove the
``towards the asymptotic direction'' condition of Theorems 1.2 and 1.3 of
\cite{CR}.

Let $\tau^{\hat v}_1$ be the first renewal time in the direction $\hat
v$, its precise definition is recalled in the next section. We prove
the following tail estimate on renewal times, which improves
Proposition 5.1 of \cite{CR}.

\begin{theo}
\label{tails-improved}
Let $l \in\mathbb S^{d-1}$, $\beta> 0$ and $M \geq15d+5$. Assume
that $(P)_M\mid l$ is satisfied and that $(E')_\beta$ holds (cf. (\ref
{kappa}), (\ref{finiteness})).
Then
\[
\limsup_{u \to\infty} (\log u)^{-1} \log P_0
\bigl( \tau_1 ^{\hat v} > u \bigr) \leq-\beta.
\]
\end{theo}

The\vspace*{1pt} condition $(E'_\beta)$ is sharp in a sense that is made precise in
Remark \ref{sharpness} below. Together with previous results of
Sznitman, Zerner, Sepp\"al\"ainen and Rassoul-Agha, cf. \cite
{SZ,Z,Sz00,RAS}, it implies the following.

\begin{theo}[(Law of large numbers)]\label{LLN-improved}
Consider a random walk in an i.i.d. environment in dimensions $d \geq
2$. Let $l \in\mathbb S^{d-1}$ and $M \geq15 d+5$. Assume that the
random walk satisfies condition $(P)_M\mid l$ and the ellipticity condition
$(E')_1$.
Then the random walk is ballistic in direction $l$ and there is a $v
\in\mathbb{R}^d$, $v \neq0$ such that
\[
\lim_{n \to\infty} \frac{X_n}{n} = v, \qquad P_0\mbox{-a.s.}
\]
\end{theo}

\begin{theo}[(Central limit theorems)]\label{CLT-improved}
Consider a random walk in an i.i.d. environment in dimensions $d \geq
2$. Let $ l \in\mathbb S^{d-1} $ and $ M \geq15 d+5 $. Assume that
the random walk satisfies condition $(P)_M\mid l$.
\begin{longlist}[(a)]
\item[(a)] (Annealed central limit theorem.) If $(E')_2$ is satisfied, then
\[
\varepsilon^{1/2} \bigl(X _{ [\varepsilon^{-1} n] } - \bigl[
\varepsilon^{-1} n\bigr] v\bigr)
\]
converges in law under $P_0$ as $\varepsilon\to0$ to a Brownian
motion with non-degenerate covariance matrix.
\item[(b)] (Quenched central limit theorem.) If $(E')_{176d}$ is
satisfied, then $\mathbb{P}$-a.s. we have that
\[
\varepsilon^{1/2} \bigl(X _{ [\varepsilon^{-1} n] } - \bigl[
\varepsilon^{-1} n\bigr] v\bigr)
\]
converges in law under $P_{0,\omega}$ as $\varepsilon\to0$ to a
Brownian motion with non-degenerate covariance matrix.
\end{longlist}
\end{theo}

Removing the ``towards the asymptotic direction'' is a real improvement:
in Section~\ref{example-dirichlet}, we will give some examples of
environments (in the class of Dirichlet environments) that satisfy
$(E'_\beta)$ but not towards the asymptotic direction. For those
environments, our new theorems allow to prove a LLN or CLT.
Furthermore, our final goal would be to get a ballisticity condition
that depends only locally on the environment (i.e., a condition that
depends only on the law of the environment at one point). Condition
$(E')_\beta$ is local, whereas $(E')_\beta$ towards the asymptotic
dimension is not: removing the ``towards the asymptotic direction'' is
then a first step in this direction. Ideally, we would also need to get
rid of condition $(P)_M\mid l$, that is not local either. This is a much
more difficult problem, not solved even in the uniformly elliptic case.

\begin{rem}\label{sharpness}
The condition of Theorem \ref{tails-improved} is sharp under the
following assumption on the tail behavior of the environment at one
site: there
exists some $(\beta_{e})_{e\in U}$, $\beta_e\geq0$, and a positive
constant $C>1$ such that for all $e\in U$
\begin{equation}
\label{sharpness-condition} C^{-1} \biggl( \prod_{e'\in U, e'\neq e}
t_{e'}^{\beta_{e'}} \biggr) \leq\mathbb{P} \bigl( \omega
\bigl(0,e'\bigr)\leq t_{e'}, \forall e'\in
U, e'\neq e \bigr) \leq C \biggl( \prod
_{e'\in U, e'\neq e} t_{e'}^{\beta_{e'}} \biggr)
\end{equation}
for all $(t_{e'})_{e'\in U\setminus\{e\}}$, $ 0 \leq t_{e'} \le1 $.
Indeed, in this case we easily see that $(E')_\beta$ is satisfied if
and only if
$\beta< 2 \sum_{e'} \beta_{e'} - \max_{e\in U} ( \beta_e +
\beta_{-e} )$.
On the other hand, if
$ \beta\geq2 \sum_{e'} \beta_{e'} - \max_{e\in U} ( \beta_e
+ \beta_{-e} )$
then $\mathbb{E} ( (\tau_1^{\hat v})^\beta) = \infty$.
Indeed, consider a direction $e_0$ which realizes the maximum in $\max
_{e\in U} ( \beta_e + \beta_{-e} )$ and set $K=\{0,e_0\}
$. We denote by $\partial_+ K$ the set of edges that exit the set $K$
composed of the edges $\{(0,e)\}_{e\neq e_0}$ and $\{(e_0, e)\}_{e\neq-e_0}$.
For small $t>0$, under the condition that $\omega(x,y-x)\leq t$ for
all $(x,y)\in\partial_+ K$ we have $P_{0,\omega} (T_K\geq n)\geq
(1-(2d-1)t)^n$. Hence,
\begin{eqnarray*}
P_0(T_K \geq n) %
&\geq& \bigl(1-(2d-1)/n
\bigr)^n \mathbb{P} \bigl( \omega(x,y-x) \leq{1/ n}, \forall(x,y)\in
\partial_+K \bigr)
\\
&\geq&\bigl(1-(2d-1)/n\bigr)^n C^{-1} n^{- ( \sum_{e'} \beta
_{e'}  - ( \beta_{e_0} + \beta_{-e_0} ) )}
\end{eqnarray*}
which implies that $E_0(T_K^\beta)=\infty$. Since $T_K$ is clearly a
lower bound for the first renewal time it gives the result.

Dirichlet environment (cf. the next section) is a typical example of
environment that satisfies condition (\ref{sharpness-condition}).
\end{rem}

\begin{rem}
\label{PM}
Theorem 1.1 of \cite{CR} states that for i.i.d. environments in
dimensions $d\geq2$ satisfying the ellipticity condition $(E')_0$, the
polynomial condition $(P)_M\mid l$ (for $ l \in\mathbb S^{d-1} $ and $M
\geq15 d+5$) is equivalent to Sznitman's condition $(T')\mid l$ (see,
e.g., \cite{Sz01} for the definition).
We can therefore replace $(P)_M\mid l$ by $(T')\mid l$ in the
statements of
Theorems \ref{LLN-improved} and \ref{CLT-improved}.
\end{rem}

\subsection{New examples of random walks satisfying the polynomial condition}
In this article, we also introduce new examples of RWRE in environments
which are not uniformly elliptic and which satisfy the polynomial
condition $(P)_M$ for $M\geq15d+5$. In Section~\ref{example1}, we
prove the polynomial condition for a subset of marginal nestling random
walks, including a particular two-dimensional environment introduced by
Campos and Ram\'\i rez in \cite{CR}. In Section~\ref
{example-dirichlet}, we prove the polynomial condition for a class of
random walks in Dirichlet random environments which do not necessarily
satisfy Kalikow's condition. In both cases, we present the case of
environments for which our new Theorems~\ref{LLN-improved} and \ref
{CLT-improved} prove necessary to study the behaviour of the walks.

\subsubsection{Example within the class of marginal nestling random walks}
\label{example1}

Following Sznitman \cite{Sz00}, we say that a law $\mathbb{P}$ on
$\Omega$ is \textit{marginal nestling} if the convex hull $K_o$ of the
support of the law of
\[
d(0,\omega):=\sum_{e\in U}\omega(0,e)e
\]
is such that $0\in\partial K_o$. We will prove in Section~\ref
{last-section1} that a certain subset of the marginal nestling laws
satisfies the polynomial condition.

\begin{theo}
\label{marginal}
Consider an elliptic law $\mathbb{P}$ under which $\{\omega(x)\dvt x\in
\mathbb Z^d\}$ are i.i.d. Assume that there exists an $r>1$ such that
$\omega(0,e_1)=r\omega(0,e_{1+d})$. Then the polynomial condition
$(P)_M\mid e_1$ is satisfied for some $M\geq15d+5$.
\end{theo}

\begin{rem}
This theorem is valid for all i.i.d. elliptic environments satisfying $
\omega(0,e_1) = r \omega(0,e_{1+d}) $, including uniformly elliptic
environments. However, the environments are marginal nestling only in
the non-uniformly elliptic case.
\end{rem}

The above result includes an example suggested in \cite{CR}, by Campos
and Ram\'{\i}rez, of an environment which satisfies the polynomial
condition and for which the random walk is directionally transient but
not ballistic. They showed that on this environment, $ (E') _{\alpha}
$ is satisfied for $\alpha$ smaller but arbitrarily close to $1$, and
that the walk is transient but not ballistic in a given direction. The
proof that this environment satisfies the polynomial condition was left
for a future work.

Let us define the environment introduced in \cite{CR}. Let $ \varphi$
be any random variable taking values on the interval $ (0, 1/4) $ and
such that the expected value of $ \varphi^{-1/2} $ is infinite, while
for every $ \varepsilon> 0 $, the expected value of $ \varphi^{ -
(1/2 - \varepsilon) } $ is finite. Let $X$ be a Bernoulli random
variable of parameter $1/2$. We now define $ \omega(0,e_1) = 2 \varphi
$, $ \omega(0,-e_1) = \varphi$, $ \omega(0,e_2)= X \varphi+ (1-X)
(1 - 4 \varphi) $ and $ \omega(0,-e_2) = X (1 - 4 \varphi) + (1-X)
\varphi$.

For every $\varepsilon>0$, this environment satisfies
$(E')_{1-\varepsilon}$: traps can appear because the random walk can
get caught on two edges of the type $ ( x, e_2 ), ( x + e_2, -e_2 )
$. Furthermore, it is transient in direction $e_1$ but not ballistic in
that direction.

\subsubsection{Examples within the class of Dirichlet random environments}\label{example-dirichlet}

Random Walks in Dirichlet Environment (RWDE) are interesting because of
the analytical simplifications they offer, and because of their link
with reinforced random walks. Indeed, the annealed law of a RWDE
corresponds to the law of a linearly directed-edge reinforced random
walk \cite{ES,P}.

Given a family of positive weights $ (\beta_1, \dots, \beta_{2d}) $,
a random i.i.d. Dirichlet environment is a law on $\Omega$ constructed
by choosing independently at each site $x \in\mathbb{Z}^d$ the values
of $ ( \omega( x, e_i ) ) _{i \in[\![1,2d]\!]} $
according to a Dirichlet law with parameters $ (\beta_1, \dots, \beta
_{2d})$. That is, at each site we choose independently a law with density
\[
\frac{ \Gamma( \sum_{i=1} ^{2d} \beta_i ) }{ \prod_{i=1} ^{2d} \Gamma
( \beta_i ) } \Biggl( \prod_{i=1}
^{2d} x_i ^{\beta_i - 1} \Biggr) \,\mathrm{d}x_1 \cdots \,\mathrm{d}x_{2d-1}
\]
on the simplex $\{(x_1,\dots, x_{2d}) \in\,]0,1]^{2d}, \sum_{i=1}
^{2d} x_i = 1 \}$. Here $\Gamma$ denotes the Gamma function $ \Gamma
(\beta) = \int_0 ^\infty t ^{\beta- 1} \mathrm{e}^{-t} \,\mathrm{d}t$, and $\,\mathrm{d}x_1 \cdots
\,\mathrm{d}x_{2d-1}$ represents the image of the Lebesgue measure on $ \mathbb
{R}^ {2d-1} $ by the application $( x_1, \dots,x_{2d-1} ) \to( x_1,
\dots, x_{2d-1}, 1- x_1 - \cdots- x_{2d-1} ) $. Obviously, the law
does not depend on the specific role of $ x_{2d}$.

\begin{rem}
Given a Dirichlet law of parameters $ ( \beta_1, \dots, \beta_{2d} )
$, the ellipticity condition $ (E') _\beta$ is satisfied if and only if
\[
\kappa\bigl( ( \beta_1, \dots, \beta_{2d} ) \bigr) = 2
\Biggl( \sum_{i=1} ^{2d}
\beta_i \Biggr) - \max_{i=1,\dots, d} (\beta_i
+ \beta_{i+d}) > \beta.
\]
As stated in Remark~\ref{sharpness}, this ellipticity condition is
optimal to get Theorem~\ref{tails-improved} in the case of Dirichlet
environments. Remark that for Dirichlet environments, for all $\beta>
0 $, $ (E') _\beta$ is much sharper that $ (E') _\beta$ towards the
asymptotic direction.
Indeed, the result of Tournier \cite{T2} gives us the explicit value
of $\hat{v}$ in the case of Dirichlet laws: $ \hat{v} = \frac{\sum
_{i=1} ^{2d} \beta_i e_i}{\llVert \sum_{i=1} ^{2d} \beta_i e_i
\rrVert }$.
Without loss of generality, we can assume that $\beta_i\ge\beta
_{i+d}$ for
$1\le i\le d$. This implies that $e_i\cdot\hat v\ge0$ for $1\le i\le d$.
If we define $\tilde\beta_i:=\min_{1\le j\le d}\beta_j$
and $\tilde\beta_{i+d}:=\min(\tilde\beta_i,\beta_{i+d})$
for $1\le i\le d$, we can see that
$(E')_\beta$ is satisfied towards the asymptotic direction
if and only if
$ \kappa( \{ \tilde\beta_i\dvt1\le i\le2d \}
) > \beta$.
\end{rem}

In the case of RWDE, it has been proved that Kalikow's condition, and
thus the $ (T') $ condition, is satisfied whenever
\begin{equation}
\label{kalikow} \sum_{i=1} ^d \llvert
\beta_i - \beta_{i+d} \rrvert> 1
\end{equation}
(see Enriquez and Sabot in \cite{ES2} and Tournier in \cite{T}). The
characterization of Kalikow's condition in terms of the parameters of a
RWDE remains an open question.
On the other hand, we believe that for RWDE condition $(T')$ is
satisfied if and only if
$ \max_{1 \leq i \leq d} \llvert \beta_i - \beta_{i+d} \rrvert > 0 $.
Nevertheless, in this article we are able to prove the following result.

\begin{theo}
\label{theorem-dirichlet}
Let $\beta_1, \beta_2,\ldots,\beta_d,\beta_{d+2},\ldots,\beta
_{2d}$ be fixed positive numbers. Then, there exists an $\varepsilon
\in(0,1)$ depending on these numbers such that if $\beta_{1+d}$ is
chosen so that $ \beta_{1+d} \leq\varepsilon$, the Random Walk in
Dirichlet Environment with parameters $ ( \beta_1, \dots, \beta_{2d}
) $ satisfies condition $(P)_M\mid e_1$ for $ M \geq15d + 5 $.
\end{theo}

Theorem \ref{theorem-dirichlet} gives as a corollary new examples of
RWDE which are ballistic in dimension $d=2$ since they do not
correspond to ranges of the parameters satisfying condition (\ref
{kalikow}) of Tournier \cite{T} and Sabot and Enriquez \cite{ES2}
(see the following remark for the case $d\geq3$).
Indeed, by Theorem \ref{LLN-improved}, if
\[
2 \sum_{i=1}^{2d}\beta_i-\max
_{1\leq i\leq d}(\beta_i+\beta_{i+d})>1
\]
and one of the parameters $\{\beta_i\dvt1\leq i\leq d\}$ is small enough,
the walk is ballistic.

\begin{rem}
In dimension $d\geq3$, in \cite{S,Bouchet}, precise conditions on
the existence of an invariant measure viewed from the particle
absolutely continuous with respect to the law have been given; this
allows to characterize completely the parameters for which there is
ballisticity, but it fails to give information on the $(T')$ condition
and on the tails of renewal times. It also fails to give a CLT.

Theorem~\ref{CLT-improved} then gives us annealed CLTs for Dirichlet
laws when the parameters $ ( \beta_1, \dots, \beta_{2d} ) $ satisfy
$2 \sum_{i=1}^{2d}\beta_i-\max_{1\leq i\leq d}(\beta_i+\beta
_{i+d})>2$ along with condition (\ref{kalikow}) or the hypothesis of
Theorem~\ref{theorem-dirichlet}.
\end{rem}

\begin{rem}
For the Dirichlet laws in dimension $d=2$ with parameters $ ( \beta_1,
\dots, \beta_{4} ) $ satisfying $2 \sum_{i=1}^{4}\beta_i-\max_{1\leq
i\leq2}(\beta_i+\beta_{i+d})>1$, with one of the parameters
$\{\beta_i\dvt1\leq i\leq4\}$ small enough, but for which $(E')_1$ is
not satisfied toward the asymptotic direction, our Theorem~\ref
{LLN-improved} gives the ballisticity when the results of \cite{CR}
would not have been enough.

For the Dirichlet laws in dimension $d\geq2$ with parameters $ ( \beta
_1, \dots, \beta_{2d} ) $ satisfying $ 2 \sum_{i=1} ^{2d} \beta_i -
\max_{1\leq i\leq d}( \beta_i + \beta_{i+d} ) > 2 $, with condition
(\ref{kalikow}) or the hypothesis of Theorem~\ref{theorem-dirichlet},
but for which $(E')_2$ is not satisfied toward the asymptotic
direction, our Theorem~\ref{CLT-improved} gives the annealed CLT when
the results of \cite{CR} would not have been enough.

This illustrates the relevance of having removed the ``toward the
asymptotic direction'' hypothesis in Theorem~\ref{tails-improved}.
\end{rem}

\section{First tools for the proofs}
In this section, we will introduce some tools that will prove necessary
for the proof of Theorem~\ref{tails-improved}.

\subsection{Regeneration times}

The proofs in \cite{CR} are based on finding bounds on the
regeneration times. We thus begin by giving the definition and some
results about the regeneration times with respect to a fixed direction
$l$. In the following, we suppose that the walk is transient in
direction $l$.

We define $\{ \theta_n\dvt n \geq1 \}$ as the canonical time shift on
$( \mathbb{Z}^d )^\mathbb{N}$. For $ l \in\mathbb S^{d-1} $ and $u
\geq0$, we define the time
\[
T ^l _u:= \min\{ n \geq0\dvt X_n \cdot l \geq
u \}.
\]

Set
\begin{equation}
\label{defa} a > 2 \sqrt{d}
\end{equation}
and
\[
D^l:= \min\{ n \geq0\dvt X_n \cdot l < X_0
\cdot l \}.
\]

We define
\begin{eqnarray*}
S_0&:=& 0, \qquad M_0:= X_0 \cdot l,
\\
S_1&:=& T^l _{M_0+a}, \qquad R_1:=
D^l \circ\theta_{S_1} + S_1,
\\
M_1&:=& \max\{ X_n \cdot l\dvt 0 \leq n \leq
R_1 \},
\end{eqnarray*}
and recursively for $k \geq1$,
\begin{eqnarray*}
S_{k+1}&:=& T ^l _{M_k+a}, \qquad
R_{k+1}:= D^l \circ\theta_{S_{k+1}} +
S_{k+1},
\\
M_{k+1}&:=& \max\{ X_n \cdot l\dvt 0 \leq n \leq
R_{k+1} \}.
\end{eqnarray*}

The first regeneration time is then defined as
\[
\tau_1:= \min\{ k \geq1\dvt S_k <\infty,
R_k = \infty\}.
\]
We can now define recursively in $n$ the $(n + 1)$th regeneration time
$ \tau_{n+1} $ as $ \tau_1 (X_\cdot) + \tau_n ( X_{\tau
_1+\cdot} - X _{\tau_1} ) $. We will occasionally write $ \tau
_1^l, \tau_2^l, \ldots$ to emphasize the dependence on the chosen direction.

\begin{rem}
The condition (\ref{defa}) on $a$ is only necessary to prove the
non-degeneracy of the covariance matrix of part $(a)$ of Theorem \ref
{CLT-improved}.
\end{rem}

It is a standard fact (see, e.g., Sznitman and Zerner \cite
{SZ}) to show that under the assumption of transience in direction $l$,
the sequence $ ( ( \tau_1, X _{(\tau_1+\cdot) \land\tau_2} - X
_{\tau_1} ), ( \tau_2 - \tau_1, X _{(\tau_2+\cdot) \land\tau_3}
- X _{\tau_2}), \ldots) $ is independent and (except for its first
term) i.i.d. Its law is the same as the law of $\tau_1$ with respect
to the conditional probability measure $ P_0 ( \cdot\mid D^l = \infty)$.

Those regeneration times are particularly useful to us because of the
two following theorems.

\begin{theo}[(Sznitman and Zerner \cite{SZ}, Zerner \cite{Z}, Sznitman
\cite{Sz00})]\label{LLN-CLT-ren}
Consider a RWRE in an elliptic i.i.d. environment. Let $l \in\mathbb
S^{d-1}$ and assume that there is
a neighbourhood $V$ of $l$ such that for every $l' \in V$ the random
walk is transient in the direction $l'$. Then there is a deterministic
$v$ such that $P_0$-a.s.
\[
\lim_{n \to\infty} \frac{X_n}{n} = v.
\]
Furthermore, the following are satisfied.
\begin{longlist}[(a)]
\item[(a)] If $ E_0 ( \tau_1 ) < \infty$, the walk is ballistic and
$v \neq0$.
\item[(b)] If $ E_0 ( \tau_1^2 ) < \infty$,
\[
\varepsilon^{1/2} \bigl( X_{ [\varepsilon^{-1} n] } - \bigl[\varepsilon
^{-1} n\bigr] v \bigr)
\]
converges in law under $P_0$ to a Brownian motion with non-degenerate
covariance matrix.
\end{longlist}
\end{theo}

\begin{theo}[(Rassoul-Agha and Sepp\"{a}l\"{a}inen \cite{RAS})]\label{CLTq-ren}
Consider a RWRE in an elliptic i.i.d. environment. Take $ l \in\mathbb
S^{d-1} $ and let $\tau_1$ be the corresponding regeneration time.
Assume that
\[
E_0 \bigl( \tau_1 ^p \bigr)< \infty,
\]
for some $ p > 176 d $. Then $\mathbb{P}$-a.s. we have that
\[
\varepsilon^{1/2} \bigl( X_{ [\varepsilon^{-1} n] } - \bigl[\varepsilon
^{-1} n\bigr] v \bigr)
\]
converges in law under $P_{0,\omega}$ to a Brownian motion with
non-degenerate covariance matrix.
\end{theo}

\subsection{Atypical quenched exit estimate}

The proof of Theorem \ref{tails-improved} is based on an atypical
quenched exit estimate proved in \cite{CR}. We will also need this
result, and thus recall it in this section. Let us first introduce some
notations.

Without loss of generality, we can assume that $e_1$ is contained in
the open half-space defined by the asymptotic direction so that
\[
\hat v\cdot e_1>0.
\]
We define the hyperplane:
\[
H:=\bigl\{ x \in\mathbb{R}^d\dvt x \cdot e_1 = 0 \bigr\}.
\]
Let $ P:= P_{\hat v} $ be the projection on the asymptotic direction
along the hyperplane $H$ defined for $ z \in\mathbb{Z}^d $ by
\[
P(z):= \biggl( \frac{ z \cdot e_1 }{ \hat{v} \cdot e_1 } \biggr) \hat{v},
\]
and $ Q:= Q_l $ be the projection of $z$ on $H$ along $\hat v$ so that
\[
Q(z)\dvt = z - P(z).
\]

Now, for $ x \in\mathbb{Z}^d$, $ \beta> 0 $, $\rho>0$ and $ L > 0
$, we define the tilted boxes with
respect to the asymptotic direction $\hat v$ by:
\begin{equation}
\label{tilted-box} B_{\beta,L}(x):= \bigl\{ y \in\mathbb{Z}^d
\mbox{ s.t. } - L^{\beta} < (y-x) \cdot e_1 < L
\mbox{ and }\bigl\llVert Q(y-x) \bigr\rrVert_\infty< \rho
L^{\beta} \bigr\}
\end{equation}
and their front boundary by
\[
\partial^{+} B_{\beta,L} (x):= \bigl\{ y \in\partial
B_{\beta,L}(x)\mbox{ s.t. } (y -x) \cdot e_1 = L \bigr\}.
\]

We have the following.

\begin{prop}[(Atypical Quenched Exit Estimate, Proposition 4.1 of \cite{CR})]
\label{aqee}
Assume there exists $ \alpha>0 $ such that $ \eta_\alpha:= \max_{e\in
U} \mathbb{E} ( ( \frac{1}{\omega(0,e)}
)^\alpha) < \infty$. Take $ M \geq15d+5 $ such that $(P)_M\mid l$
is satisfied. Let $\beta_0 \in(1/2,1)$, $\beta\in( \frac
{\beta_0+1}{2}, 1 )$ and $ \zeta\in(0,\beta_0) $.
Then, for each $\gamma> 0$ we have that
\[
\limsup_{L\to\infty} L ^{-g(\beta_0,\beta,\zeta)} \log\mathbb{P} \bigl(
P_{0,\omega} \bigl( X _{ T_{B_{\beta,L}(0)} } \in\partial^{+}
B_{\beta,L}(0) \bigr) \leq \mathrm{e}^{-\gamma L^\beta} \bigr) < 0,
\]
where
\[
g (\beta_0,\beta,\zeta):= \min\bigl\{ \beta+ \zeta, 3\beta- 2 +
(d-1) (\beta-\beta_0) \bigr\}.
\]
\end{prop}

\subsection{Some results on flows}

The main tool that enables us to improve the results of \cite{CR} is
the use of flows and max-flow min-cut theorems. We need some
definitions and properties that we will detail in this section. In the
following, we consider a finite directed graph $ G = ( V, E ) $, where
$V$ is the set of vertices and $E$ is the set of edges. For all $ e \in
E $, we denote by $\underline{e}$ and $\overline{e}$ the vertices
that are the tail and head of the edge $e$ (the edge $e$ goes from
$\underline{e}$ to $\overline{e}$).

\begin{de}
We consider a finite directed graph $ G = ( V, E ) $. A flow from a
set $ A \subset V $ to a set $ Z \subset V $ is a non-negative function
$ \theta\dvtx E \to\mathbb{R}_+ $ such that:
\begin{itemize}
\item$ \forall x \in(A \cup Z) ^c $, $\operatorname{div} \theta(x)
= 0 $.
\item$ \forall x \in A $, $ \operatorname{div} \theta(x) \geq0 $.
\item$ \forall x \in Z $, $ \operatorname{div} \theta(x) \leq0 $.
\end{itemize}
Where the divergence operator is $ \operatorname{div}\dvtx \mathbb{R}^E
\to\mathbb{R}^V $ such that for all $ x \in V $,
\[
\operatorname{div} \theta(x) = \sum_{e \in E,\underline{e}=x} \theta(e)
- \sum_{e \in E,\overline{e}=x} \theta(e).
\]

A unit flow from $A$ to $Z$ is a flow such that $ \sum_{x \in A}
\operatorname{div} \theta(x) = 1 $. (Then we have also $ \sum_{x \in
Z} \operatorname{div} \theta(x) = -1 $.)
\end{de}

We will need the following generalized version of the max-flow min-cut theorem.

\begin{prop}[(Proposition 1 of \cite{S})]\label{max-flow-min-cut}
Let $ G = (V,E) $ be a finite directed graph. Let $ ( c(e) )_{e \in E}
$ be a set of non-negative reals (called capacities). Let $x_0$ be a
vertex and $ ( p_x ) _{x \in V} $ be a set of non-negative reals. There
exists a non-negative function $ \theta\dvtx E \to\mathbb{R}_+ $ such that
\begin{equation}
\label{divergence-condition} \operatorname{div} \theta= \sum_{ x \in V }
p_x ( \delta_{x_0} - \delta_x )
\end{equation}
and
\begin{equation}
\label{capacities-condition} \forall e \in E,\qquad \theta(e) \leq c(e)
\end{equation}
if and only if for all subset $ K \subset V $ containing $x_0$ we have
\begin{equation}
\label{cutset-big-enough} c ( \partial_+ K ) \geq\sum_{ x \in K ^c }
p_x,
\end{equation}
where $ \partial_+ K = \{ e \in E, \underline{e} \in K, \overline{e}
\in K^c \} $ and $ c ( \partial_+ K ) = \sum_{ e \in\partial_+ K }
c(e) $. The same is true if we restrict the condition (\ref
{cutset-big-enough}) to the subsets $K$ such that any $ y \in K $ can
be reached from $x_0$ following a directed path in $K$.
\end{prop}

We will give here an idea of the proof, that explains why we call this
result a generalized version of the classical max-flow min-cut theorem.
The complete proof can be found in \cite{S}.

\begin{pf*}{Idea of the proof}
If $\theta$ satisfies (\ref{divergence-condition}) and (\ref
{capacities-condition}), then
\[
\sum_{e, \underline{e} \in K, \overline{e} \in K^c} \theta(e) - \sum
_{e, \overline{e} \in K, \underline{e} \in K^c} \theta(e) = \sum_{x \in K}
\operatorname{div} \theta(x) = \sum_{x \in K^c}
p_x.
\]
It implies (\ref{cutset-big-enough}) by (\ref{capacities-condition})
and positivity of $\theta$.

The reversed implication is an easy consequence of the classical
max-flow min-cut theorem on finite directed graphs (see, e.g.,
\cite{LP} Section~3.1). If $ ( c(e) )_{e \in E} $ satisfies (\ref
{cutset-big-enough}), we consider the new graph $\tilde{G} = ( V \cup
\{ \delta\}, \tilde{E} )$, where
\[
\tilde{E} = E \cup\bigl\{ (x,\delta), x \in V \bigr\}.
\]
We define a new set of capacities $( \tilde{c} (e) ) _{e \in\tilde
{E}}$ where $ c(e) = \tilde{c} (e) $ for $ e \in E $ and $ \tilde{c}
( ( x, \delta) ) = p_x $. The strategy is to apply the max-flow
min-cut theorem with capacities $\tilde{c}$ and with source~$x_0$ and
sink~$\delta$.
It gives a flow $ \tilde{\theta} $ on $ \tilde{G} $ between $ x_0 $
and $ \delta$ with strength $ \sum_{x \in V} p_x $ and such that $
\tilde{\theta} \leq\tilde{c} $. The function $\theta$ obtained by
restriction of $ \tilde{\theta} $ to $ E $ satisfies (\ref
{capacities-condition}) and (\ref{divergence-condition}).
\end{pf*}

For the proof of Theorem \ref{tails-improved}, we will consider the
oriented graph $ ( \mathbb{Z}^d, E _{\mathbb{Z}^d} ) $ where $ E
_{\mathbb{Z}^d}:= \{ (x,y) \in(\mathbb{Z}^d)^2$ s.t. $
\llvert x - y \rrvert _1 = 1 \} $. This graph is not finite, but we
will only
consider flows with compact support ($ \theta(e) = 0 $ for all $e$
except in a finite subset of $E _{\mathbb{Z}^d}$). We can then proceed
as if the graph were finite, and use the previous definition and proposition.

\section{Proof of Theorem \texorpdfstring{\protect\ref{tails-improved}}{1}}

We will prove Theorem \ref{tails-improved} using the atypical quenched
exit estimate Proposition \ref{aqee}. Let us give a rough idea of the
proof. We first show that the event $\{\tau_1>u\}$ is concentrated on
the event that the random walk does not exit a box of side $(C\log
u)^{1/\beta}$, for an appropriate choice of $C$, before time
$u$. Now on this last event, necessarily, the walk must visit some
point of this box at least $N_u:=u/(C\log u)^{{d}/{\beta}}$ times.
But due to Proposition \ref{aqee} and the strong Markov property, the
probability that this point is visited $N_u$ times is less\vspace*{1pt} than $
(1-\frac{1}{u^{1-\varepsilon}} )^{u/(C\log u)^{{d}/{\beta
}}}$, for some $\varepsilon>0$ which depends on the choice of $C$.
This last quantity tends quickly to $0$, and then the dominant term
bounding $P_0(\tau_1>u)$ will be the probability to exit the box of
side $(C\log u)^{1/\beta}$ before time~$u$.

Let $l \in\mathbb S^{d-1}$, $\beta> 0$ and $M \geq15d+5$. Assume
that $(P)_M\mid l$ is satisfied and that $(E')_\beta$ holds.

Let us take a rotation $\hat R$ such that $ \hat R(e_1) = \hat v $. We
fix $ \beta' \in( \frac{5}{6}, 1 )$, $M>0$ and for
simplicity we will write $\tau_1$ instead of $\tau_1 ^{\hat{v}}$.

For $ u > 1 $, take
\begin{eqnarray*}
L &=& L(u):= \biggl( \frac{1}{4M\sqrt{d}} \biggr) ^{{1}/{\beta
'}} ( \log u)
^{{1}/{\beta'}},
\\
C_L&:=& \hat{R} \biggl( \biggl[ - \frac{L}{2(\hat{v}\cdot e_1)},
\frac{L}{2(\hat{v}\cdot e_1)} \biggr] ^d \biggr) \cap\mathbb{Z}^d.
\end{eqnarray*}

Following the proof of Proposition 5.1 in \cite{CR}, we write
\[
P_0 (\tau_1 >u) \leq P_0 (
\tau_1 > u, T_{C_{L(u)}} \leq\tau_1 ) + \mathbb{E}
\bigl( F_1 ^c, P_{0,\omega} (T_{C_{L(u)}}>u )
\bigr) + \mathbb{P}(F_1),
\]
with
\[
F_1:= \biggl\{ \omega\in\Omega\dvt t_{\omega} (
C_{L(u)} ) > \frac{u}{(\log u)^{{1}/{\beta'}}} \biggr\}
\]
and
\[
t_{\omega}(A):= \min\biggl\{ n \geq0\dvt \sup_{x}
P_{x,\omega} ( T_A > n ) \leq\frac{1}{2} \biggr\}.
\]

As in \cite{CR}, the term $P_0 ( \tau_1 > u, T_{C_{L(u)}} \leq
\tau_1 )$ is bounded thanks to condition $ (P) _M \mid l $, and
the term $\mathbb{E} ( F_1 ^c, P_{0,\omega}
(T_{C_{L(u)}}>u ) )$ is bounded thanks to the strong
Markov property. This part of the original proof is not modified, so we
will not give more details here. It gives the existence for every $
\gamma\in( \beta', 1 ) $ of a constant $c > 0$ such that:
\[
P_0 (\tau_1 >u) \leq\frac{ \mathrm{e}^{- c L (u) ^\gamma} }{c} + \biggl(
\frac{1}{2} \biggr) ^{ \lfloor( \log u) ^{{1}/{\beta'}}
\rfloor} + \mathbb{P}(F_1).
\]
It only remains to show that we can find a constant $ C > 0$ such that
$ \mathbb{P}(F_1) \leq C u ^{-\beta} $ for $u$ big enough.

For each $\omega\in\Omega$, still as in \cite{CR}, there exists
$x_0 \in C_{L(u)}$ such that
\[
P _{x_0,\omega} ( \tilde{H}_{x_0} > T_{C_{L(u)}} ) \leq
\frac{ 2
\llvert C_{L(u)}\rrvert }{ t_{\omega}(C_{L(u)}) },
\]
where for $y \in\mathbb{Z}^d$, $\tilde{H}_y = \min\{ n \geq1\dvtx X_n
=y \}$. It gives
\[
\mathbb{P}(F_1) \leq\mathbb{P} \biggl( \omega\in\Omega\mbox
{ s.t. } \exists x_0 \in C_{L(u)}\mbox{ s.t. }
P_{x_0,\omega} ( \tilde{H}_{x_0} > T _{C_{L(u)}} ) \leq
\frac{2 (\log u) ^{{1}/{\beta'}}}{u} \llvert C_{L(u)} \rrvert\biggr).
\]

We define for each point $x \in C_{L(u)}$ a point $y_x$, closest from $
x + 2 \frac{L^{\beta'}}{\hat{v}\cdot e_1} \hat{v}$. To bound $
\mathbb{P}(F_1)$, we will need paths that go from $x$ to $y_x$ with
probability big enough and the atypical quenched exit estimate
(Proposition \ref{aqee}).

Define:
\[
N:= \frac{ \llvert \hat{v} \rrvert \log u }{ 2 M \sqrt{d}
(\hat{v}\cdot e_1) }.
\]
It is straightforward that
\[
N - 1 \leq\llvert y_x - x \rrvert_1 \leq N+1.
\]

The following of the proof will be developed in three parts: first, we
will construct unit flows $ \theta_{i,x} $ going from $ \{ x, x + e_i
\} $ to $ \{ y_x, y_x + e_i \} $, for all $x \in C_{L(u)}$. Then we
will construct paths with those flows, and use the atypical quenched
exit estimate to bound $\mathbb{P}(F_1)$ in the case that those paths
are big enough. We will conclude by bounding the probability that the
paths are not big enough.

\subsection{Construction of the flows \texorpdfstring{$\theta_{i,x}$}{theta{i,x}}}
We consider the oriented graph $ ( \mathbb{Z}^d, E _{\mathbb{Z}^d} )
$ where $ E _{\mathbb{Z}^d}:= \{ (x,y) \in(\mathbb{Z}^d)^2$
s.t. $\llvert x - y \rrvert _1 = 1 \} $. We want to construct
unit flows $
\theta_{i,x} $ going from $ \{ x, x + e_i \} $ to $ \{ y_x, y_x +
e_i \} $, for all $x \in C_{L(u)}$. But there are additional
constraints, as we will need them to construct paths that have a
probability big enough. The aim of this section is to prove the
following proposition.

\begin{prop}\label{existence-of-flows}
For all $x \in C_{L(u)}$, for all $ \alpha_1, \dots, \alpha_{2d} $
positive constants, there exists $2d$ unit flows $ \theta_{i,x}\dvt E
_{\mathbb{Z}^d} \to\mathbb{R}_+ $, respectively, going from $ \{ x,
x + e_i \} $ to $ \{ y_x, y_x + e_i \} $, such that:
\begin{equation}
\label{bound-on-theta} \forall e \in E _{\mathbb{Z}^d},\qquad \theta_{i,x}(e)
\leq
\frac
{\alpha(e)}{\kappa_i},
\end{equation}
where $ \kappa_i:= 2 \sum_{j=1} ^{2d} \alpha_j - (\alpha_i +
\alpha_{i+d})$, and $ \alpha(e):= \alpha_j $ for $e$ of the type $(z,e_j)$.

Furthermore, we can construct $\theta_{i,x}$ with a finite support,
and in a way that allows to find $\gamma$ and $ S \subset E _{\mathbb
{Z}^d} $, $\llvert S \rrvert $ independent of $u$, such that
$\theta_{i,x}(e) \kappa_i \leq\gamma< \alpha(e) $ for all $ e \in
S ^c $.
\end{prop}

We will construct the $ \theta_{i,x} $ to prove their existences. For
this, we need three steps. Let $ B(x,R) $ be the box of $ \mathbb{Z}^d
$ of center $x$ and radius $R$, and $ B_i(x,R) $ be the same box, where
the vertices $x$ and $x+e_i$ are merged (and we suppress the edge
between them). We note $ E _{B(x,R)}:= \{ (x,y) \in E _{\mathbb{Z}^d}
\cap( B(x,R) )^2 \} $ and $ E _{B_i(x,R)}:= \{ (x,y) \in E _{\mathbb
{Z}^d} \cap( B_i(x,R) )^2 \} $ the corresponding sets of edges. We
will construct a unit flow in the graph $ ( B_i(x,R), E _{B_i(x,R)} )
$ from $ \{ x, x + e_i \} $ to $ B_i(x,R) ^c $, a unit flow in the
graph $ ( B_i(y_x,R), E _{B_i(y_x,R)} ) $ from $ B_i(y_x,R) ^c $ to $
\{ y_x, y_x + e_i \} $, and then connect them. At each step, we will
ensure that condition (\ref{bound-on-theta}) is fulfilled.

\begin{longlist}[\textit{First step}.]
\item[\textit{First step}.] Construction of a unit flow from $ \{ x, x +
e_i \} $ to $ B_i(x,R) ^c $:

\begin{lemma}\label{flow-in-the-ball}
Set $x \in C_{L(u)}$, and $ \alpha_1, \dots, \alpha_{2d} $ positive
constants. If $ R \geq\frac{ \max_{i} \kappa_i }{ \min_j \alpha
_j }$, there exists $2d$ unit flows $ \theta_{i,x}\dvtx E _{B_i(x,R)}
\to
\mathbb{R}_+ $ such that:
\[
\operatorname{div} \theta_{i,x} = \sum_{z \in\partial B_i(x,R)}
\frac{1}{\llvert \partial B_i(x,R) \rrvert } ( \delta_x - \delta_z)
\]
and
\[
\forall e \in E _{B_i(x,R)},\qquad \theta_{i,x}(e) \leq
\frac{\alpha
(e)}{\kappa_i},
\]
where $ \partial B_i(x,R) = \{ z \in B_i(x,R)$ that has a neighbour in
$B_i(x,R) ^c \} $.
\end{lemma}

The divergence condition ensures that the flow will be a unit flow,
that it goes from $x$, and that it leaves $B_i(x,R)$ uniformly on the
boundary of the box.

\begin{pf*}{Proof of Lemma \ref{flow-in-the-ball}}
The result is a simple application of Proposition \ref
{max-flow-min-cut}. We fix $ x \in C_{L(u)} $ and $i$ between $1$ and
$2d$. Define $ p_z = \frac{1}{\llvert \partial B_i(x,R) \rrvert } $
if $ z \in\partial B_i (x,R) $, $ p_z = 0 $ if $ z \notin
\partial B_i (x,R) $.

To prove the result, we only have to check that $ \forall K \subset
B_i(x,R) $ containing $x$, $ \sum_{e \in\partial_+ K} \frac{\alpha
(e)}{\kappa_i} \geq\sum_{z \notin K} p_z $, where $ \partial_+ K =
\{ e \in E _{B_i(x,R)}$ s.t. $\underline{e} \in K$ and $\overline{e}
\notin K \} $.

We have two cases to examine:
\begin{itemize}
\item If $ K \cap\partial B_i (x,R) = \varnothing$, $\sum_{z \notin
K} p_z = 1$. We then need $ \sum_{e \in\partial_+ K} \alpha(e) \geq
\kappa_i $. For $ K = \{ x \} $, $ \sum_{e \in\partial_+ K} \alpha
(e) = \kappa_i $ as we merged $x$ and $x+e_i$. For bigger $K$, we
consider for all $j \neq i $ the paths $ (x + n e_j) _{n \in\mathbb
{N}} $ and for all $j \neq i + d $ the paths $ (x + e_i + n e_j) _{n
\in\mathbb{N}} $. They intersect the boundary of $K$ in $2d+1$
different points, and the exit directions give us the corresponding
$\alpha_j$, that sum to $ \kappa_i $. It gives that $ \sum_{e \in
\partial_+ K} \alpha(e) \geq\kappa_i $.
\item If $ K \cap\partial B_i (x,R) \neq\varnothing$, $\sum_{z
\notin K} p_z < 1$. As $K$ contains\vspace*{2pt} a path from $x$ to $ \partial B_i
(x,R) $, $\sum_{e \in\partial_+ K} \frac{\alpha(e)}{\kappa_i}
\geq\frac{R \min_j (\alpha_j + \alpha_{j+d})}{\kappa_i}$. It is
bigger than $1$ thanks to the hypothesis on $R$. It gives the result.\quad\qed
\end{itemize}\noqed
\end{pf*}

\item[\textit{Second step}.] By the same way, we construct a flow $ \theta
_{i,x}\dvtx E _{B_i(y_x,R)} \to\mathbb{R}_+ $ such that
\[
\operatorname{div} \theta_{i,x} = \sum_{z \in\partial B_i(y_x,R)}
\frac{1}{\llvert \partial B_i(y_x,R) \rrvert } ( \delta_z - \delta_{y_x})
\]
and
\[
\forall e \in E _{B_i(y_x,R)},\qquad \theta_{i,x}(e) \leq
\frac{\alpha
(e)}{\kappa_i}.
\]

\item[\textit{Third step}.] We will join the flows on $ E _{B_i(x,R)} $ and
$ E _{B_i(y_x,R)} $ with simple paths, to get a flow on $ E _{\mathbb
{Z}^d} $. Take $ R \geq\frac{ \max_{i} \kappa_i }{ \min_j \alpha
_j }$, and make sure that $ \frac{1}{\llvert \partial B(x,R)
\rrvert } < \frac{\alpha(e)}{\kappa_i} $ for all $ e \in E
_{\mathbb{Z}^d} $ (always possible by taking $R$ big enough, $R$
depends only on the $\alpha_i$ and the dimension).

We can find $\llvert \partial B(x,R)\rrvert $ simple paths $ \pi_j
\subset E
_{\mathbb{Z}^d} $ satisfying:
\begin{itemize}
\item$\forall j$, $ \pi_j $ connects a point of $\partial B(x,R)$ to
a point of $\partial B(y_x,R)$.
\item$\forall j$, $ \pi_j $ stays outside of $B(x,R)$ and $B(y_x,R)$,
except from the departure and arrival points.
\item If two paths intersect, they perform jumps in different direction
after the intersection (no edge is used by two paths). If $ (x, e_i) $
is in a path, then $(x+e_i, -e_i )$ is not in any path.
\item The number of steps of each path is close to $N$: there exists
constants $K_1$ and $K_2$ independent of $u$ such that the length of $
\pi_j $ is smaller than $K_1 N + K_2$.
\end{itemize}
(E.g., we can use the paths $ \pi^{(i,j)} $ page~45 of \cite
{CR}, and make them exit the ball $B(x,R)$ instead of $ \{ x, x + e_i
\} $).

For all $i$, $\partial B_i(x,R) = \partial B(x,R)$ and $\partial
B_i(y_x,R) = \partial B(y_x,R)$ as soon as $ R > 1 $. By construction,
$ - \operatorname{div} \theta_{i,x} ( z_1 ) = \operatorname{div}
\theta_{i,x} ( z_2 ) = \frac{1}{\llvert \partial B(x,R) \rrvert } $
for any $ z_1 \in\partial B_i(x,R)$ and $ z_2 \in\partial
B_i(y_x,R) $. We can then join the flows of the first two steps by
defining a flow $ \theta_{i,x} (e) = \frac{1}{\llvert \partial
B(y_x,R) \rrvert } $ for all $e \in\pi_j$ (and $0$ on all the
other edges of $ E _{\mathbb{Z}^d} $).

We have thus constructed a unit flow $ \theta_{i,x} $ on $ E _{\mathbb
{Z}^d} $, from $ \{ x, x + e_i \} $ to $ \{ y_x, y_x + e_i \} $,
satisfying (\ref{bound-on-theta}) ((\ref{bound-on-theta}) is
satisfied on $ E _{B_i(x,R)} $ and $ E _{B_i(y_x,R)} $ as $ R \geq
\frac{ \max_{i} \kappa_i }{ \min_j \alpha_j } $ thanks to Lemma
\ref{flow-in-the-ball}, and outside those balls as $ \frac{1}{\llvert
\partial B(x,R) \rrvert } < \frac{\alpha(e)}{\kappa_i} $
for all $ e \in E _{\mathbb{Z}^d} $). It concludes the proof of the
first part of Proposition \ref{existence-of-flows}.

As $ \theta_{i,x} (e) = 0 $ out of the finite set $ E _{B_i(x,R)} \cup
E _{B_i(y_x,R)} \cup\{ e \in\pi_j, 1 \leq j \leq\llvert \partial
B(x,R)\rrvert
\} $, the flow has a finite support. And as we made sure that $ \frac
{1}{\llvert \partial B(x,R) \rrvert } < \frac{\alpha
(e)}{\kappa_i} $, we can take $S = B(x,R) \cup B(y_x,R)$ and $ \gamma
= \frac{\kappa_i}{\llvert \partial B(x,R) \rrvert } $ to
conclude the proof.
\end{longlist}

\subsection{Bounds for $\mathbb{P}(F_1)$}

We apply Proposition \ref{existence-of-flows} for the $ \alpha_1,
\dots, \alpha_{2d} $ of the definition of $(E')_\beta$ (see (\ref
{kappa}) and (\ref{finiteness})). It gives flows $ \theta_{i,x} $ on
$ E _{\mathbb{Z}^d} $, constructed as in the previous section.

We can decompose a given $ \theta_{i,x} $ (for $i$ and $x$ fixed) in a
finite set of weighted paths, each path starting from $x$ or $x+e_i$
and arriving to $y_x$ or $y_x + e_i$. It suffices to choose a path
$\sigma$ where the flow is always positive, to give it a weight $
p_\sigma:= \min_{e \in\sigma} \theta_{i,x} (e) > 0 $ and to
iterate with the new flow $ \theta(e):= \theta_{i,x} (e) - p _\sigma
\mathbh{1} _{ e \in\sigma} $.

The weight $p_\sigma$ of a path $\sigma$ then satisfies: for all $e
\in E _{\mathbb{Z}^d} $, $ \theta_{i,x} (e) = \sum_{\sigma~\mathrm
{containing}~e} p_\sigma$. As $\theta_{i,x}$ is a unit flow, we get $
\sum_{\sigma~\mathrm{path~of}~\theta_{i,x}} p_\sigma= 1 $. We will\vspace*{1pt}
use those weights in the next section, to prove that those paths are
``big enough'' with high probability.

We now introduce:
\[
F_{2,i} = \biggl\{ \omega\in\Omega\mbox{ s.t. } \forall x \in
C_{L(u)}, \exists\sigma\mbox{ path of } \theta_{i,x},
\omega_{\sigma}:= \prod_{e \in\sigma}
\omega_e > u ^{({1}/{M}) - 1} \biggr\},
\]
where we recall that $M$ is the parameter for condition $(P)_M $, and
\[
F_{2} = \bigcap_{i=1} ^{2d}
F_{2,i}.
\]

Define
\[
F_3:= F_2 \cap\biggl\{ \omega\in\Omega\mbox{ s.t.
} \exists x_0 \in C_{L(u)}\mbox{ s.t. }
P_{x_0,\omega} ( \tilde{H}_{x_0} > T _{C_{L(u)}} ) \leq
\frac{2 (\log u) ^{{1}/{\beta'}}}{u} \llvert C_{L(u)} \rrvert\biggr
\}.
\]

We get immediately:
\[
\mathbb{P}(F_1) \leq\mathbb{P}(F_3) + \mathbb{P}
\bigl(F_2 ^c\bigr).
\]

It gives two new terms to bound. We start by bounding $\mathbb
{P}(F_3)$. For this we will use the same method as in \cite{CR}: on
the event $F_3$, for all $ 1 \leq i \leq2d $ we can use a path $\sigma
$ of $ \theta_{i,x} $ to join $x$ or $x+e_i$ to $y_x$ or $y_x + e_i$.
It gives:
\begin{eqnarray*}
\omega( x_0, e_i ) u ^{ (1/M) - 1 } \min
_{z \in\{ y_{x_0},
y_{x_0} + e_i \} } P_{z,\omega} ( T _{C_{L(u)}} <
H_{x_0} ) &\leq& P_{x_0,\omega} ( T _{C_{L(u)}} <
\tilde{H}_{x_0} )
\\
&\leq& \frac{2
(\log u) ^{{1}/{\beta'}}}{u} \llvert C_{L(u)}
\rrvert,
\end{eqnarray*}
where the factor $\omega( x_0, e_i )$ corresponds to the probability
of jumping from $x$ to $x+e_i$, in the case where the path $\sigma$
starts from $x+e_i$.\vspace*{1pt}

As $ \sum_{i=1} ^{2d} \omega( x_0, e_i ) = 1 $, it gives
\[
u ^{ (1/M) - 1 } \min_{z \in V( y_{x_0} ) } P_{z,\omega} ( T
_{C_{L(u)}} < H_{x_0} ) \leq\frac{4d (\log u) ^{{1}/{\beta
'}}}{u} \llvert
C_{L(u)} \rrvert,
\]
where $ V( y_{x_0} ):= \{ y_{x_0}, (y_{x_0} + e_i) _{i=1, \dots, 2d}
\} $.

In particular, on $F_3$, we can see that for $u$ large enough $ V(
y_{x_0} ) \subset C_{L(u)}$. As a result, on $F_3$, we have for $u$
large enough
\begin{eqnarray*}
\min_{z \in V( y_{x_0} ) } P_{z,\omega} ( X _{T_{z} + U _{\beta',L}}
\cdot
e_1 > z \cdot e_1 ) &\leq& \min_{z \in V( y_{x_0} ) }
P_{z,\omega} ( T _{C_{L(u)}} < H_{x_0} )
\\
&\leq&
\frac{1}{ u ^{1/(2M)}} = \mathrm{e}^{ -2 \sqrt d L(u) ^{\beta'} },
\end{eqnarray*}
where
\[
U _{\beta',L}:= \bigl\{ x \in\mathbb{Z}^d\dvt -L ^{\beta'}
< x \cdot e_1 < L \bigr\}.
\]
From this and using the translation invariance of the measure $\mathbb
{P}$, we conclude that:
\begin{eqnarray*}
&&  \mathbb{P} \biggl( \exists x_0 \in C_{L(u)}
\mbox{ s.t. } P_{x_0,\omega} ( \tilde{H}_{x_0} > T
_{C_{L(u)}} ) \leq\frac{4d
(\log u) ^{{1}/{\beta'}}}{u} \llvert C_{L(u)} \rrvert,
F_2 \biggr)
\\
&&\quad \leq\mathbb{P} \Bigl( \exists x_0 \in
C_{L(u)}\mbox{ s.t. } \min_{z \in V( y_{x_0} ) }
P_{z,\omega} ( X _{T_{z + U _{\beta',L}}} \cdot e_1 > z \cdot
e_1 ) \leq \mathrm{e}^{ -2 \sqrt d L(u) ^{\beta'} } \Bigr)
\\
&&\quad \leq(2d+1) \llvert C_{L(u)} \rrvert\mathbb{P} \bigl(
P_{0,\omega} ( X _{ T _{U _{\beta', L(u)}} } \cdot e_1 > 0 ) \leq \mathrm{e}^{ -2 \sqrt d L(u) ^{\beta'} } \bigr)
\\
&&\quad \leq(2d+1) \llvert C_{L(u)} \rrvert\mathbb{P} \bigl(
P_{0,\omega} ( X _{ T _{B _{\beta', L(u)}} } \cdot e_1 > 0 ) \leq \mathrm{e}^{ -2 \sqrt d L(u) ^{\beta'} } \bigr),
\end{eqnarray*}
where the tilted box $ B _{\beta', L(u)} $ is defined as in (\ref
{tilted-box}).

We conclude with the atypical quenched exit estimate (Proposition \ref
{aqee}): there exists a constant $c>0$ such that for each $ \beta_0
\in(\frac{1}{2},1)$ one has:
\[
\mathbb{P}(F_3) \leq\frac{1}{c} \mathrm{e}^{-c L(u) ^{g(\beta_0, \beta',
\zeta)}},
\]
where $g(\beta_0,\beta',\zeta)$ is defined as in Proposition \ref{aqee}.

Note that for each $ \beta' \in( \frac{5}{6},1 )$ there
exists a $ \beta_0 \in( \frac{1}{2},\beta) $ such that
for every $ \zeta\in( 0,\frac{1}{2} )$ one has $ g (
\beta_0, \beta', \zeta) > \beta' $. Therefore, replacing $L$ by its
value, we proved that there exists $c > 0$ such that:
\[
\mathbb{P}(F_3) \leq c u^{-\beta}.
\]

\subsection{Bound for $\mathbb{P}(F_2^c)$}

To conclude the bound for $ \mathbb{P}(F_1) $ and the proof of Theorem
\ref{tails-improved}, it only remains to control $\mathbb{P}(F_2
^c)$. It is in this section that we will use the conditions that were
imposed on $ \theta_{i,x} $ during the construction of the flows, as
well as condition $ (E')_\beta$
\begin{eqnarray*}
\mathbb{P}\bigl(F_2 ^c\bigr)
&\leq& \sum_{i=1} ^{2d} \mathbb{P}
\bigl(F_{2,i}^c\bigr)
\\
& \leq& \sum_{i=1} ^{2d} \sum
_{x \in C_{L(u)}} \mathbb{P}\bigl(\forall\sigma
\mbox{ path of } \theta_{i,x}, \omega_{\sigma} \leq u
^{(1/M) - 1} \bigr).  %
\end{eqnarray*}

As $\theta_{i,x}$ is a unit flow, if $\forall\sigma$ path of $\theta
_{i,x}$, $\omega_{\sigma} \leq u ^{(1/M) - 1}$ then:
\[
\sum_{\sigma~\mathrm{path~of}~\theta_{i,x}} p_\sigma\omega_{\sigma}
\leq u ^{(1/M) - 1} \sum_{\sigma~\mathrm
{path~of}~\theta_{i,x}} p_\sigma=
u ^{(1/M) - 1}.
\]
Jensen's inequality then gives:
\[
\prod_{e \in E _{\mathbb{Z}^d} } \omega_{e}^{\theta_{i,x}(e)} =
\prod_{\sigma~\mathrm{path~of}~\theta_{i,x}} \omega_{\sigma} ^{p_\sigma}
\leq u ^{(1/M) - 1}.
\]
It allows to write:
\begin{eqnarray*}
\mathbb{P}\bigl(F_2 ^c\bigr) &\leq& \sum_{i=1} ^{2d} \sum
_{x \in C_{L(u)}} \mathbb{P} \biggl( \prod_{e \in E _{\mathbb{Z}^d} }
\omega_{e}^{\theta_{i,x}(e)} \leq u ^{(1/M) - 1} \biggr)
\\
&\leq& \sum_{i=1} ^{2d} \sum
_{x \in C_{L(u)}} \frac{ \mathbb{E}
( \prod_{e \in E _{\mathbb{Z}^d} } \omega_{e}^{-\kappa_i \theta
_{i,x}(e)} ) }{ u ^{-\kappa_i (({1}/{M}) - 1)} }.
\end{eqnarray*}

We will use the integrability given by the flows to bound the
expectations. The independence of the environment gives (for $i$ and
$x$ fixed):
\begin{eqnarray*}
\mathbb{E} \biggl( \prod_e
\omega_{e}^{-\kappa_i \theta_{i,x}(e)} \biggr)
&=& \prod
_{z \in\mathbb{Z}^d} \mathbb{E} \biggl( \prod_{e~\mathrm
{s.t.}~\underline{e} = z}
\omega_{e}^{-\kappa_i \theta_{i,x}(e)} \biggr)
\\
&=& \prod_{z \in S} \mathbb{E} \biggl( \prod
_{e~\mathrm{s.t.}~\underline{e} = z} \omega_{e}^{-\kappa_i \theta
_{i,x}(e)} \biggr)
\prod_{z \notin S} \mathbb{E} \biggl( \prod
_{e~\mathrm{s.t.}~\underline{e} = z} \omega_{e}^{-\kappa_i \theta
_{i,x}(e)} \biggr),
\end{eqnarray*}
where we recall that $S = B(x,R) \cup B(y_x,R)$.

As $ \theta_{i,x} $ satisfies (\ref{bound-on-theta}), the ellipticity
condition $ (E') _\beta$ gives that each of the expectations $\mathbb
{E} ( \prod_{e~\mathrm{s.t.}~\underline{e} = z} \omega_{e}^{-\kappa_i \theta_{i,x}(e)} )$ are finite.

By construction $ \llvert S \rrvert $ is finite and does not
depend on $u$: $\prod_{z \in S} \mathbb{E} ( \prod_{e~\mathrm
{s.t.}~\underline{e} = z} \omega_{e}^{-\kappa_i \theta_{i,x}(e)}
)$ is a finite constant independent on $u$.

It remains to deal with the case of $z \notin S$. As we chose $R$ to
get $ \theta_{i,x}(e) \kappa_i < \gamma$ for the edges outside $S$,
and thanks to the bounds on the number of edges with positive flow
(there is a finite number of paths, and each path has a bounded
length), we have:
\[
\prod_{z \notin S} \mathbb{E} \biggl( \prod
_{e~\mathrm{s.t.}~\underline{e} = z} \omega_{e}^{-\kappa_i \theta
_{i,x}(e)} \biggr) \leq
\mathbb{E} \biggl( \prod_{e~\mathrm{s.t.}~\underline{e} = 0} \omega_{e}^{- \gamma} \biggr) ^{c_1 N + c_2},
\]
where $c_1$ and $c_2$ are positive constants, independent of $u$. Then,
putting all of those bounds together,
\begin{eqnarray*}
\mathbb{P}\bigl(F_2 ^c\bigr)
&\leq& \sum_{i=1} ^{2d} \sum
_{x \in C_{L(u)}} C_1 C_2 ^{C_3 N} u
^{\kappa_i ((1/M) - 1)}
\leq \sum_{i=1} ^{2d} \sum
_{x \in C_{L(u)}} C_4 u ^{ ((C_5 +\kappa_i)/{M}) - \kappa_i}
\\
&\leq& C_6 (\log u)^{C_7} u ^{(C_8/M) - \min_i \kappa_i},
\end{eqnarray*}
where all the constants $C_i$ are positive and do not depend on $u$. As
Remark \ref{PM} tells us that we can choose $M$ as large as we want,
we can get $\frac{C_8}{M}$ as small as we want.

Then we can find a constant $C > 0$ such that $\mathbb{P}(F_2 ^c) \leq
C u ^{-\beta} $ for $u$ big enough. It concludes the proof.

\section{New examples of random walks satisfying the polynomial condition}
\label{last-section1}
\subsection{Proof of Theorem \texorpdfstring{\protect\ref{marginal}}{4}}
Consider the box $B_{e_1,L,\tilde L}$ for
$ \tilde{L} = 70 L^3 $. We want to find some integer
$L > c_0 $ such that
\[
P _0 ( X_{ T_{ B_{e_1,L,\tilde L} } } \cdot e_1 < L ) \leq
\frac{1}{L^M},
\]
for some $M\geq15d+5$.
We first decompose this probability according to whether the exit point
of the random walk from the box $B_{e_1,L,\tilde L}$ is on the bottom
or on one of the sides of the box, so that,
\begin{eqnarray*}
&& P _0 ( X_{ T_{ B_{e_1,L,\tilde L} } } \cdot e_1 < L )
\\
&&\quad = P _0 ( X_{ T_{ B_{e_1,L,\tilde L} } } \cdot e_1 = - L
) + \sum_{i=2} ^d \bigl( P
_0 ( X_{ T_{ B_{e_1,L,\tilde L} } } \cdot e_i = \tilde{L} ) + P
_0 ( X_{ T_{ B_{e_1,L,\tilde
L} } } \cdot e_i = - \tilde{L} ) \bigr).
\end{eqnarray*}
We will first bound the probability to exit through the sides. We do
the computations for $ P _0 ( X_{ T_{ B_{e_1,L,\tilde L} } }
\cdot e_2 = \tilde{L} ) $ but the\vspace*{-1pt} other term can be dealt with
in the same way. Suppose that $ X_{ T_{ B_{e_1,L,\tilde L} } } \cdot
e_2 = \tilde{L} $, and define $ n_0, \dots, n_{ \tilde{L} - 1 } $
the finite hitting times of new levels in direction $e_2$ as follows:
\[
n _k:= \min\{ n \geq0~\mathrm{s.t.}~X_n \cdot
e_2 \geq k \}.
\]
To simplify notation define $\varphi(x):=\omega(x,e_{1+d})$.
We now choose a constant $ 1 > \delta> 0 $, and we will call ``good point''
any $ x \in\mathbb{Z}^2 $ such that $ \varphi( x ) > \delta$. We
define $p:= \mathbb{P}( \varphi( x ) > \delta)$. Note that
$p$ does not depend on $x$ since the
environment is i.i.d., so that it depends only on $ \delta$ and the
law of~$ \varphi$.

We now introduce the event that a great number of the $X _{n _k}$ are
good points:
\[
C_1:= \biggl\{ X_{ T_{ B_{e_1,L,\tilde L} } } \cdot e_2 = \tilde{L}
\mbox{ and at least }\frac{p}{2} \tilde{L}\mbox{ of the } X
_{n _k}, 1 \leq k \leq\tilde{L} - 1, \mbox{ are good} \biggr\}.
\]
We get immediately
\[
P _0 ( X_{ T_{ B_{e_1,L,\tilde L} } } \cdot e_2 = \tilde{L} ) = P
_0 ( C _1 ) + P _0 \bigl( \{
X_{ T_{
B_{e_1,L,\tilde L} } } \cdot e_2 = \tilde{L} \} \cap\bigl(C_1
^c\bigr) \bigr).
\]
By construction of the $ X _{n _k} $ and independence of the
environment, and with $Z$ an independent random variable following a
binomial law of parameters $ p $ and $ \tilde{L} $, we can bound the
second term of the sum:
\begin{eqnarray*}
P _0 \bigl( \{ X_{ T_{ B_{e_1,L,\tilde L} } } \cdot e_2 =
\tilde{L} \} \cap\bigl(C_1 ^c\bigr) \bigr) &\leq&  P \biggl( Z \leq\frac{p}{2} \tilde{L} \biggr)
\\
& \leq& \exp\biggl( - 2 \frac{ ( p \tilde{L} - p \tilde{L} / 2 ) ^2
}{ \tilde{L} } \biggr)
\\
& =& \exp\biggl( - \frac{ p ^2 \tilde{L} }{ 2 } \biggr),
\end{eqnarray*}
where the last inequality is Hoeffding's inequality.

It only remains to bound $ P _0 ( C_1 ) $. For that, we introduce the
following new event
\[
C_2:= C_1 \cap\biggl\{ X _{n _k + 1} - X
_{n _k} = e_1\mbox{ for at least } \frac{\delta p}{4}
\tilde{L}\mbox{ of the good } X _{n _k} \biggr\},
\]
that states that the walk goes often in direction $e_1$ just after
reaching a $X _{n _k}$ that is a good point.

We can then write
\[
P _0 ( C_1 ) = P _0 ( C_2 ) +
P _0 \bigl( C_1 \cap\bigl(C_2
^c\bigr) \bigr).
\]
To bound $ P _0 ( C_1 \cap(C_2 ^c) ) $, we use the uniform bound
``$ \varphi( x ) > \delta$'' for good points that gives us that
$ \omega( x, e_1 ) > r \delta$ on those points. Getting $Z'$ an
independent random variable following a binomial law of parameters $ r
\delta$ and $ \lfloor\frac{p}{2} \tilde{L} \rfloor$,
it gives:
\begin{eqnarray*}
P _0 \bigl( C_1 \cap\bigl(C_2
^c\bigr) \bigr) %
& \leq& P \biggl( Z' \leq
\frac{\delta p}{4} \tilde{L} \biggr)
\leq \exp\biggl( - p \delta^2 \tilde{L} \biggl( r -
\frac{1}{2} \biggr)^2 \biggr).
\end{eqnarray*}

It only remains to bound $ P _0 ( C_2 ) $. Set $n^+$ (resp.,
$n^-$) the total number of jumps in direction $e_1$ (resp.,
$-e_1$) before exiting the box $B_{e_1,L,\tilde L}$. We will need a
third new event
\[
C_3:= \biggl\{ n^+ \geq\frac{1+r}{2} n^- \biggr\},
\]
that allows us to write
\[
P _0 ( C_2 ) = P _0 ( C_2 \cap
C_3 ) + P _0 \bigl( C_2 \cap
\bigl(C_3 ^c\bigr)\bigr).
\]
First, notice that for $L$ big enough, $ C_2 \cap C_3 = \varnothing$. Indeed,
$C_1$ implies that we exit the box $B_{e_1,L,\tilde L}$ by the side ``$
x \cdot
e_2 = \tilde{L} $''. Now, since the vertical displacement
of the walk before exiting the box $B_{e_1,L,\tilde L}$ is $n^+-n^-$,
on the event $C_3$ we know that this displacement is at
least equal to $\frac{r-1}{r+1}n^+$. Therefore, since on $C_2$ the
walk makes at least
$ \frac{\delta p}{4} \tilde{L} = \frac{35 \delta p}{2} L ^3 $
moves in the direction $ e_1$, on $C_2\cap C_3$ its vertical
displacement before exiting the box is at least
$\frac{35\delta p(r-1)}{2(r+1)} L^3$.
Since on $C_2\cap C_3$ the walk exits the box by the ``$ x \cdot e_2 =
\tilde{L} $'' side,
we see that for $L$ larger than $L_1:= \sqrt{ \frac{2(1+r)}{35
\delta p (r-1)} }$ the event $C_2\cap C_3$ is empty.

We now want to bound $ P _0 ( C_2 \cap(C_3 ^c)) $
\begin{eqnarray*}
P _0 \bigl( C_2 \cap\bigl(C_3
^c\bigr)\bigr) %
&\leq& P _0 \biggl( n^+ \geq
\frac{\delta p}{4} \tilde{L}\mbox{ and } n^+ < \frac{1+r}{2} n^-
\biggr)
\\
&\leq& P _0 \biggl( n^+ + n^- \geq\frac{\delta p}{4}
\tilde{L}\mbox{ and }\bigl(n^+ + n^-\bigr) \frac{2}{3+r} < n^-
\biggr).
\end{eqnarray*}
Now note that whenever we go through a vertical edge from a point $x$,
the law of the environment tells us that it is an edge $ (x,e_1) $ with
probability $\frac{r}{1+r} $, and $ (x,-e_1) $ with probability $\frac
{1}{1+r} $. Then, defining $Z''$ as a random variable following a
binomial law of parameters $ \frac{1}{1+r} $ and $ \lfloor\frac
{p \delta}{4} \tilde{L} \rfloor$, we have the bound:
\begin{eqnarray*}
P _0 \bigl( C_2 \cap\bigl(C_3
^c\bigr)\bigr) %
&\leq& P \biggl( Z''
\geq\frac{2}{3+r} \frac{p \delta}{4} \tilde{L} \biggr)
\\
&\leq& \exp\biggl( - 2 p \delta\tilde{L} \biggl( \frac{1}{3+r} -
\frac{2}{p \delta\tilde{L}(1+r)} \biggr) ^2 \biggr),
\end{eqnarray*}
where\vspace*{1pt} we need $ 1 \leq\frac{ (r+1) p \delta\tilde{L}}{2(3+r)} $ to
apply Hoeffding's inequality in the last inequality. We can find $ L_2
$ such that this is true for $ L \geq L_2 $.

Choose $ M \geq15d+5$. By putting all of our previous bounds together,
we finally get, for all $ L \geq L_2 $,
\begin{eqnarray*}
&& P _0 ( X_{ T_{ B_{e_1,L,\tilde L} } } \cdot e_2 =
\tilde{L} )
\\
&&\quad \leq\exp\biggl( - \frac{ p ^2 \tilde{L} }{ 2 } \biggr) + \exp\biggl( - p
\delta^2 \tilde{L} \biggl( r - \frac{1}{2}
\biggr)^2 \biggr) + \exp\biggl( - 2 p \delta\tilde{L} \biggl(
\frac{1}{3+r} - \frac{2}{p \delta\tilde{L}(1+r)} \biggr) ^2 \biggr),
\end{eqnarray*}
where we recall that $ \tilde{L} = 70 L^3 $, $ \delta> 0 $ and $ p =
\mathbb{P}(
\varphi(x) > \delta) $. Then, for any choice of $ \delta$, we can
find $ L_3 \geq\max( c_0, L_1, L_2 ) $ such that for all $ L \geq
L_3 $,
\[
\sum_{i=2} ^d \bigl( P _0 (
X_{ T_{ B_{e_1,L,\tilde L} } } \cdot e_i = \tilde{L} ) + P _0 (
X_{ T_{ B_{e_1,L,\tilde
L} } } \cdot e_i = -\tilde{L} ) \bigr) \leq
\frac{1}{2 L
^M}.
\]

We now only need to bound $ P _0 ( X_{ T_{ B_{e_1,L,\tilde L} } }
\cdot e_1 =
-L )$ to prove $ (P)_M \mid e_1 $. We will use again the notations $n^+$
(resp., $n^-$) for the total number of jumps in direction $e_1$
(resp., $-e_1$) before exiting the box $B_{e_1,L,\tilde L}$. Suppose
that $ X_{ T_{ B_{e_1,L,\tilde L} } } \cdot e_1 = -L $. Then
necessarily $ n^+
< n^- $, which gives $ n^+ < \frac{ n^+ + n^- }{2} $. As $n^+$
conditioned to
$n^++n^-$ follows a binomial law of parameters $ \frac{r}{1+r} > \frac
{1}{2} $ and $ n^+ + n^- $, Hoeffding's inequality gives the bound:
\[
P_0 \biggl( n^+ < \frac{ n^+ + n^- }{2} \bigg| n^++n^- \biggr) \leq\exp
\biggl( - 2 \bigl( n^+ + n^- \bigr) \biggl( \frac{r}{1+r} -
\frac{1}{2} \biggr) ^2 \biggr).
\]
But $ X_{ T_{ B_{e_1,L,\tilde L} } } \cdot e_1 = -L $ also gives that
necessarily, $ n^- \geq L $. Then
\begin{eqnarray*}
P _0 ( X_{ T_{ B_{e_1,L,\tilde L} } } \cdot e_1 = -L )
&\leq& P_0 \biggl( n^+ < \frac{ n^+ + n^- }{2}
\mbox{ and }n^- \geq L \biggr)
\\
& \leq& \sum_{m=L}^\infty\exp
\biggl( - 2 m \biggl( \frac{r}{1+r} - \frac{1}{2} \biggr)
^2 \biggr).
\end{eqnarray*}
Therefore, we can find $ L_4 \geq L_3 $ such that for all $ L \geq L_4 $,
\[
P _0 ( X_{ T_{ B_{e_1,L,\tilde L} } } \cdot e_1 = -L ) \leq
\frac{1}{2 L ^M},
\]
from where we conclude that for all $ L \geq L_4 $,
\[
P _0 ( X_{ T_{ B_{e_1,L,\tilde L} } } \cdot e_1 < L ) \leq
\frac{1}{L ^M}.
\]

\subsection{Proof of Theorem \texorpdfstring{\protect\ref{theorem-dirichlet}}{5}}

It is classical to represent Dirichlet distributions with independent
gamma random variables: if $ \gamma_1, \dots, \gamma_N $ are
independent gamma random variables with parameters $ \beta_1, \dots,
\beta_N $, then $ \frac{\gamma_1}{\sum\gamma_i}, \dots, \frac
{\gamma_N}{\sum\gamma_i} $ is a Dirichlet random variable with
parameters $ ( \beta_1, \dots, \beta_N ) $. We get a restriction
property as an easy consequence of this representation (see \cite{W},
pages 179--182): if $ ( \xi_1, \dots, \xi_N ) $ is a Dirichlet
random variable with parameters $ ( \beta_1, \dots, \beta_N ) $, for
any $J$ non-empty subset of $ \{ 1, \dots, N \} $, the random variable
$ ( \frac{\xi_j}{\sum_{i \in J} \xi_i} ) _{j \in J} $
follows a Dirichlet law with parameters $ ( \beta_j ) _{j \in J} $ and
is independent of $ \sum_{i \in J} \xi_i $. This property will be
useful in the following.

We consider the box $B_{e_1,L,\tilde L}$ for $ \tilde{L} = 70 L^3 $,
and want to find some $ L > c_0 $ such that $ P _0 ( X_{ T_{
B_{e_1,L,\tilde L} } } \cdot e_1 < L ) \leq\frac{1}{L^M} $ to
prove $(P)_M \mid e_1$.

Let $ l_i:= \{ x \in\mathbb{Z}^d$ s.t. $x \cdot e_1 = i \} $ and
$t_i:= \min\{ n \geq0\dvt X_n \in l_i, X _{n+1} \notin l_i \}$. We
first consider the events that, when the walk arrives on $l_i$ for the
first time, it gets out of it by an edge in direction $e_1$ (the
alternative being getting out by an edge in direction $-e_1$):
\[
G _{1,i}:= \{ X _{ t_i +1 } - X _{ t_i } =
e_1 \}.
\]

At the point $X _{ t_i }$, we know that the walk will go either to $X
_{ t_i } + e_1$ or to $X _{ t_i } - e_1$. Thanks to the restriction
property of the Dirichlet laws, we know that $ \frac{ \omega( X _{
t_i }, e_1) }{ \omega( X _{ t_i }, e_1) + \omega( X _{ t_i }, -
e_1) } $ follows a beta law of parameters $ ( \beta_1, \beta_1 +
\beta_{1+d} ) $ and is independent from the previous trajectory of the
walk on~$l_i$. Indeed, we already know by construction that it is
independent from the environment on the other points of $l_i$, and the
restriction property gives the independence\vspace*{1pt} from $ \frac{ \omega( X
_{ t_i }, e_j) }{ 1 - \omega( X _{ t_i }, e_1) - \omega( X _{ t_i }, -
e_1) } $ for all $j \neq1, -1$, which corresponds to the law of
the previous trajectory from $X _{ t_i }$ on $l_i$. Then
\[
P _0 ( G _{1,i} ) = \frac{ \beta_1 }{ \beta_1 + \beta_{1+d} }.
\]

Now define
\[
G _1:= \bigcap_{i=1} ^{ L } G
_{1,i},
\]
and note that
\[
P _0 \bigl( G _1^c \bigr) \leq L
\frac{\beta_{1+d}}{\beta_1+\beta_{1+d}}.
\]
We can now write
\begin{eqnarray*}
&& P _0 ( X_{ T_{ B_{e_1,L,\tilde L} } } \cdot e_1 < L )
\\
&&\quad \leq P _0 \bigl( \{ X_{ T_{ B_{e_1,L,\tilde L} } } \cdot
e_1 < L \} \cap G _1 \bigr) + L\frac{\beta_{1+d}}{\beta_1+\beta_{1+d}},
\end{eqnarray*}
and we only need to bound the first term of this sum. If $G _1$ is
satisfied, the walk cannot get out of the box $B_{e_1,L,\tilde L}$ by
the ``lower boundary'' $ \{ x \in\mathbb{Z}^d$ s.t. $x \cdot e_1 = -L
\} $. Then the walk has to get out by one of the $2d-2$ ``side
boundaries'':
\[
P _0 \bigl( \{ X_{ T_{ B_{e_1,L,\tilde L} } } \cdot e_1 < L \} \cap
G _1 \bigr) = P _0 \Biggl( \bigcup
_{j=2}^d \{ X_{ T_{ B_{e_1,L,\tilde L} } } \cdot
e_j = \pm\tilde{L} \} \cap G _1 \Biggr).
\]

On the event $ \bigcup_{j=1}^d\{X_{ T_{ B_{e_1,L,\tilde L} } } \cdot
e_j =
\pm\tilde{L}\}$ define $ n_0, \dots, n_{ \tilde{L} - 1 }$ as the finite
hitting times of new levels in any direction perpendicular to $e_1$ as follows:
\[
n _k:= \min\Bigl\{ n \geq0\mbox{ s.t. } \max_{2\leq j\leq d}
\llvert X_n \cdot e_j\rrvert\geq k \Bigr\}.
\]
Let now $p = \frac{ \beta_1 }{1 + \sum_{i \neq1+d} \beta_i}$ and
consider the event
\[
G _{3}:= G_1 \cap\biggl\{ X _{n _k + 1} - X
_{n _k } = e_1\mbox{ for at least } \frac{ p}{2}
\tilde{L}\mbox{ of the points } X _{n _k} \biggr\}.
\]
Suppose\vspace*{1pt} $ \beta_{1+d} \leq1 $, then $ p \leq\mathbb{E}( \omega
(0,e_1) ) $.
Consider now a random variable $Z$ with a binomial law of
parameters $p$ and $\tilde L$. Using Hoeffding's inequality, we see that
\[
P\bigl(G_3^c\bigr)\leq P \biggl(Z\le\frac{p\tilde L}{2}
\biggr) \le\exp{ \biggl( -\frac{p^2}{2}\tilde L \biggr)}.
\]
But clearly $G_1\cap G_3=\varnothing$ for $ L \geq L_0:= \sqrt{\frac
{1}{35p}}$. Therefore, we have in this case
\[
P _0 \Biggl( \bigcup_{j=2}^d
\{ X_{ T_{ B_{e_1,L,\tilde L} } } \cdot e_j = \pm\tilde{L} \} \cap G
_1 \Biggr) \leq\exp{ \biggl( -\frac{p^2}{2}\tilde L \biggr)}.
\]
Putting the previous bounds together, we finally get for all $ L \geq
L_0 $:
\[
P _0 ( X_{ T_{ B_{e_1,L,\tilde L} } } \cdot e_1 < L )\leq L
\frac{\beta_{1+d}}{\beta_1+\beta_{1+d}}+\exp{ \biggl( -\frac
{p^2}{2}\tilde L \biggr)}.
\]
Let now $L_1$ be such that for all $L \geq L_1$
\[
\exp\biggl( -\frac{p^2}{2}\tilde L \biggr) \leq\frac{1}{2 L^M}.
\]
The constant $c_0$ (cf. (\ref{ce0})) is increasing in $\beta_{1+d}$.
Therefore, in the region where $\beta_{1+d}\le\beta_1$, it does not
depend on $\beta_{1+d}$. Call this value $c_0'$.
On the other hand, by construction, $L_0$ and $L_1$ do not depend on
$\beta_{1+d}$.
Take now $L_2:=\max\{c'_0,L_0,L_1\}$.

We can then choose $\beta_{1+d}$ (necessarily $\leq\beta_1$) so that
\[
L_2\frac{\beta_{1+d}}{\beta_1 + \beta_{1+d}}\le\frac{1}{2 L_2^M}.
\]
We then conclude that for this choice of $\beta_{1+d}$
there exists an $L\geq c_0$ such that
\[
P _0 ( X_{ T_{ B_{e_1,L_2,\tilde L} } } \cdot e_1 < L ) \leq
\frac{1}{L ^M}.
\]


\section*{Acknowledgements}
The authors would like
to thank Alexander Drewitz for pointing out the ideas
of the proof of Theorem \ref{marginal}.

This work was partially supported by Fondo Nacional de Desarrollo Cient\'\i fico y
Tecnol\'ogico grants 1100298 and 1141094, Iniciativa Cient\'\i fica Milenio NC120062 and by the ANR project MEMEMO2.


%

\printhistory
\end{document}